\documentstyle{amsppt}
\magnification=1200
\hoffset=-0.5pc
\vsize=57.2truepc
\hsize=38truepc
\nologo
\spaceskip=.5em plus.25em minus.20em
\define\batviltw{1}
\define\batvilfo{2}
\define\batavilk{3}
\define\evluwein{4}
\define\gersthtw{5}
\define\geschthr{6}
\define\getzltwo{7}
\define\hochsone{8}
\define\poiscoho{9}
\define\duality{10}
\define\husmosta{11}
\define\kosmathr{12}
\define\koszulon{13}
\define\liazutwo{14}
\define\rinehone{15}
\define\stashnin{16}
\define\weinsfte{17}
\define\xuone{18}
\define\Bobb{\Bbb}
\define\fra{\frak}
\noindent
dg-ga/9704005
\bigskip\noindent

\topmatter
\title 
Lie-Rinehart algebras, Gerstenhaber algebras,
\\
and Batalin-Vilkovisky algebras
\endtitle
\author Johannes Huebschmann
\endauthor
\affil 
Universit\'e des Sciences et Technologies de Lille
\\
UFR de Math\'ematiques
\\
F-59 655 VILLENEUVE D'ASCQ C\'edex
\\
Johannes.Huebschmann\@univ-lille1.fr
\endaffil
\abstract{For any Lie-Rinehart algebra $(A,L)$, (i) generators  for the 
Gerstenhaber algebra $\Lambda_A L$ correspond bijectively to right $(A,L)$-
connections on $A$ in such a way that (ii) B(atalin)-V(ilkovisky) structures 
correspond bijectively to right $(A,L)$-module structures on $A$. When $L$ is 
projective as an $A$-module, given an exact generator $\partial$, the homology 
of the B-V algebra $(\Lambda_A L,\partial)$ coincides with the homology of $L$ 
with coefficients in $A$ with reference to the right $(A,L)$-module structure 
determined by $\partial$. When $L$ is also of finite rank $n$, there are  
bijective correspondences between $(A,L)$-connections on $\Lambda_A^nL$ and 
right $(A,L)$-connections on $A$ and between left $(A,L)$-module structures 
on $\Lambda_A^nL$ and right $(A,L)$-module structures on $A$. Hence there are 
bijective correspondences between $(A,L)$-connections on $\Lambda_A^n L$ and 
generators for the Gerstenhaber bracket on $\Lambda_A L$ and between 
$(A,L)$-module structures on $\Lambda_A^n L$ and B-V algebra structures on 
$\Lambda_A L$. The homology of such a B-V algebra $(\Lambda_A L,\partial)$ 
coincides with the cohomology of $L$ with coefficients in $\Lambda_A^n L$, 
with reference to the left $(A,L)$-module structure determined by $\partial$.
Some applications to Poisson structures and to differential geometry are 
discussed.}
\endabstract
\subjclass{primary 17 B 55, 17 B 56, 17 B 65, 17 B 66, 
81~C~99, 81~D~07, 81 T 70
secondary 
70~H~99}
\endsubjclass
\date {April 7, 1997}
\enddate
\endtopmatter
\document
\leftheadtext{Johannes Huebschmann}
\rightheadtext{Lie-Rinehart-, Gerstenhaber-, BV-algebras}

\beginsection Introduction

It has been known for a while
that, given 
a vector space $\fra g$ over a field $k$
or, more generally,
a projective
module $\fra g$ over a commutative ring $R$ with 1,
Lie brackets on
$\fra g$
are in bijective correspondence with
differentials 
on the graded exterior  coalgebra $\Lambda'_R \fra g$
yielding a structure of a
differential graded coalgebra.
The 
homology and cohomology of $\fra g$
with coefficients in any
$\fra g$-module $M$
may then be computed from suitable complexes
involving
$\Lambda'_R \fra g$ and $M$
and a suitable additional structure,
a twisting cochain;
in particular,
the homology of
$\Lambda'_R \fra g$
is that of 
$\fra g$ with trivial coefficients.
We refrain from spelling
out details since we shall not need them;
see e.~g. \cite\husmosta\ 
for the notion of a twisting cochain and that
of a twisted object.
For a general {\it Lie-Rinehart algebra\/} $(A,L)$
(a definition will be reproduced in Section 1 below),
a description of its homology and cohomology
is no longer available in terms of twisting cochains;
in fact, 
when the action of $L$ on $A$ is non-trivial
there is no way
to endow
$\Lambda'_A L$
with a structure of a differential graded coalgebra
corresponding to the Lie-Rinehart structure.
Rather,
given $A$, a commutative algebra
over the ground ring $R$,
and an $A$-module $L$,
{\it Lie-Rinehart structures\/} on $(A,L)$
correspond bijectively
to 
{\it Gerstenhaber algebra\/}
structures on the graded exterior $A$-algebra
$\Lambda_A L$ over $L$.
In the Gerstenhaber algebra picture,
there is no differential in sight and the question
arises 
as to what
the link between
Gerstenhaber structures
and homology and cohomology of Lie-Rinehart algebras
might be.
The purpose of the present note is to show that
an
answer
and some new insight
may be obtained
by means of
the notion of a {\it Batalin-Vilkovisky algebra\/}.
This will also shed considerable light 
on these algebras themselves.
Batalin-Vilkovisky algebras
have recently become important in string theory
and elsewhere
\cite{\batviltw\ -- \batavilk},
\cite\getzltwo, \cite\liazutwo,
\cite\stashnin.
\smallskip
Let $(A,L)$ be a Lie-Rinehart algebra.
Given $A$, 
we shall also refer to $L$ as an $(R,A)$-{\it Lie algebra\/},
cf. \cite\rinehone.
Our first result, Theorem 1,
will say that {\it there is a bijective
correspondence between
right\/}
$(A,L)$-{\it connections
on\/} $A$
(in a sense made precise in Section 1 below)
{\it and\/} $R$-{\it linear operators generating the Gerstenhaber 
bracket on\/} 
$\Lambda_A L$;
moreover,
{\it under this correspondence, flat 
right\/} $(A,L)$-{\it connections
on\/} $A$,
{\it that is,
right\/} $(A,L)$-{\it module structures
on\/} $A$,
{\it correspond to operators
of square zero, that is, to differentials\/}.
\smallskip
To explain our next result,
write $U(A,L)$ for the universal algebra for
$(A,L)$
and  recall 
the complex
$(K(A,L),d)$
generalizing the usual resolution
used in 
Cartan-Chevalley-Eilenberg
Lie-algebra cohomology \cite\rinehone, \cite\poiscoho.
We refer to
$(K(A,L),d)$
as the {\it Rinehart complex\/}
for $(A,L)$;
it is an acyclic relatively free chain complex
in the category of left $U(A,L)$-modules.
When $L$ is projective as an $A$-module,
$(K(A,L),d)$ is in fact a projective resolution
of $A$ in the category of left $(A,L)$-modules.
Whether or not $L$ is projective as an $A$-module,
given  an exact generator $\partial$ for the Gerstenhaber
algebra $\Lambda_AL$, let
$A_{\partial}$ denote $A$ together
with the corresponding
right $(A,L)$-module structure
(given by 
our first result). 
Our second result, Theorem 2 below,
will say that,
{\it as a chain
complex,
the Batalin-Vilkovisky algebra\/} $(\Lambda_AL,\partial)$
{\it coincides with\/}
$(A_{\partial}\otimes_{U(A,L)}K(A,L),d)$.
Thus,
when $L$ is projective as an $A$-module,
the Batalin-Vilkovisky algebra
$(\Lambda_AL,\partial)$
computes the homology
$$
\roman H_*(L,A_\partial) \left(= \roman{Tor}_*^{U(A,L)}(A_\partial,A)\right)
$$
of $L$ with coefficients in
the right $(A,L)$-module $A_\partial$;
we mention in passing that in general
$(\Lambda_AL,\partial)$
will compute a certain relative homology
in the sense of relative homological algebra
\cite\hochsone.
\smallskip
Thereafter we consider
the special
case where, as an $A$-module,
$L$ is  projective  
of finite rank (say) $n$,
so that
$\Lambda_A^nL$ is the highest non-zero exterior power of $L$
in the category of $A$-modules.
In an earlier paper \cite\duality,
we introduced
the concept of 
{\it dualizing module\/}
for such a Lie-Rinehart algebra. 
Our next result, 
Theorem 3, will 
exploit this notion
and will
say that
{\it there is a bijective
correspondence between\/}
$(A,L)$-{\it connections on\/} $\Lambda_A^nL$
{\it and
right\/}
$(A,L)$-{\it connections
on\/} $A$;
moreover,
under this correspondence,
{\it left\/} $(A,L)$-{\it module structures on\/} $\Lambda_A^nL$,
{\it that is,  flat connections,
correspond to
right\/}
$(A,L)$-{\it module structures
on\/} $A$.
Theorem 1 and Theorem 3
together imply at once, cf. the Corollary
in Section 2, that
{\it there is a bijective correspondence between\/}
$(A,L)$-{\it connections
on\/}
$\Lambda_A^n L$
{\it and linear operators generating the Gerstenhaber 
bracket on\/}
$\Lambda_A L$;
further,
{\it under this correspondence, flat connections correspond to operators
of square zero,
that is, to differentials.\/}
This generalizes certain observations
made by Koszul \cite\koszulon\ 
and Xu \cite\xuone, cf. the remark in Section 2.
\smallskip
Finally, we combine the present results with
those
obtained in \cite\duality\
establishing
homological duality for Lie-Rinehart algebras,
in the following way:
Let
$\nabla$ be
a flat 
$(A,L)$-connection
on
$\Lambda_A^n L$,
let
$A_\nabla$
be the corresponding
right
$(A,L)$-module 
(Theorem 3),
and let
$\partial$
be the corresponding
exact generator
for the Gerstenhaber algebra $\Lambda_A L$ 
(Theorem 1).
By Theorem 2,
as a  chain complex,
the
Batalin-Vilkovisky algebra
$(\Lambda_A L,\partial)$
coincides
with
$(A_\nabla \otimes_{U(A,L)}K(A,L),d)$
and the latter computes
the homology of $L$
with coefficients in the right
$(A,L)$-module 
$A_\nabla$.
Theorem 4 below will say
that this homology and hence that of the
Batalin-Vilkovisky algebra
$(\Lambda_A L,\partial)$
{\it is naturally isomorphic to
the cohomology\/}
$$
\roman H^*(L,\Lambda_A^n L_\nabla)
\quad (=\roman{Ext}^*_{U(A,L)}(A,\Lambda_A^n L_\nabla))
$$
{\it of\/} $L$
{\it with coefficients in\/} $\Lambda_A^n L_\nabla$,
the $A$-module
$\Lambda_A^n L$
with left $(A,L)$-module structure
given by the connection $\nabla$.
\smallskip
In a sense, the approach in the present paper
relies on a careful study of the interplay between
left- and right modules
and of the resulting homological algebra
for a general Lie-Rinehart algebra or, more generally,
of the interplay between left and right connections.
Once the appropriate notions and language
have been found,
to some extent, the theory takes care of itself.
\smallskip
To some extent,
our terminology 
is the same as that in \cite\kosmathr.
I am indebted to Jim Stasheff for 
(i) a number of comments on an earlier draft which helped
improve the exposition and for
(ii) having drawn my attention
to Xu's paper \cite\xuone;
in fact, the latter prompted me to write the present paper.
In a final section we explain briefly how our results are
related to those of Xu's.
By extending our results to appropriate
strong homotopy notions
we intend to {\it tame\/} elsewhere the bracket zoo
that arose recently
in topological field theory,
cf. e.~g. \cite\stashnin.

\beginsection 1. Lie-Rinehart-, Gerstenhaber-, BV-algebras

In this section, we shall
explain the significance of 
the generator of a
Gerstenhaber structure
and that of an exact structure.
For ease of exposition we first recall the definitions.
\smallskip
Let $R$ be a commutative ring with 1
and $A$ a commutative $R$-algebra.
An $(R,A)$-{\it Lie algebra\/}
\cite\rinehone\ 
is a Lie algebra $L$ over $R$ which acts
on (the left of) $A$ 
(by derivations)
and is also an $A$-module 
satisfying suitable compatibility conditions
which generalize the usual properties
of the Lie algebra of vector fields on a smooth manifold
viewed as a module over its ring of functions;
these conditions read
$$
\align
[\alpha,a\beta]  &= \alpha(a)\beta + a [\alpha,\beta],
\\
(a\alpha)(b) &= a (\alpha(b)),
\endalign
$$ 
where
$a,b \in A$ and $\alpha, \beta \in L$.
When the emphasis is
on the pair $(A,L)$
with the mutual structure
of interaction between $A$ and $L$,
we refer to the pair $(A,L)$ as
a
{\it Lie-Rinehart} algebra;
in \cite\geschthr\ 
the terminology {\it Palais pair\/} is used.
Some relevant history
may be found the introduction and in Section 1 of \cite\poiscoho.
\smallskip
We now suppose
given an arbitrary
Lie-Rinehart algebra
$(A,L)$.
Consider the graded exterior $A$-algebra
$\Lambda_AL$
over $L$ 
where $L$ is taken concentrated in degree 1.
Write typical elements of
$\Lambda_A L$ in the form
$$
\langle \alpha_1,\alpha_2,\dots,\alpha_n \rangle,
\quad \alpha_1,\alpha_2,\dots,\alpha_n \in L.
$$
The Lie bracket $[\cdot,\cdot]$
of $L$ 
induces an $R$-linear bracket
on
$\Lambda_AL$
which endows the latter
with a structure of 
a
{\it Gerstenhaber algebra\/}.
Abusing notation somewhat, we continue to denote the resulting bracket
by
$$
[\cdot,\cdot]
\colon
\Lambda_AL
\otimes_R
\Lambda_AL
@>>>
\Lambda_AL.
$$
To get an explicit formula for it, let
$u = \langle\alpha_1, \dots, \alpha_\ell\rangle \in \Lambda^{\ell}_AL$ and
\linebreak
$v = \langle \alpha_{\ell+ 1}, \dots, \alpha_n\rangle \in \Lambda^{n-\ell}_AL$,
where $\alpha_1, \dots, \alpha_n \in L$;
then
$$
[u,v]
=
(-1)^{|u|} 
\sum_{j\leq \ell <k} (-1)^{(j+k)}
\langle\lbrack \alpha_j,\alpha_k \rbrack, 
\alpha_1,\dots \widehat{\alpha_j} \dots \widehat{\alpha_k}
\dots, \alpha_n \rangle.
\tag1.1
$$
For example when $A$ is the ring of functions
on a smooth manifold and $L$ the Lie algebra of vector fields,
$\Lambda_AL$
is the algebra of  multivector fields
and the bracket is the
{\it Schouten\/} bracket.
In general, a {\it Gerstenhaber algebra\/}
is a graded commutative algebra
$\Cal A$
together with a Lie bracket
from $\Cal A \otimes_R \Cal A$ to $\Cal A$
of degree $-1$ (in an appropriate sense:
it is a graded Lie bracket in the usual sense
when the degrees of the elements of $\Cal A$
are lowered by 1);
see \cite\geschthr\ 
where these objects are called $G$-algebras,
or \cite{\kosmathr,\,\liazutwo,\,\stashnin,\,\xuone}.
We recall from
Theorem 5 of  \cite\geschthr\ 
that (i) the assignment
to
$\Cal A$ of the pair
$(A_0, A_1)$
consisting of the homogeneous degree zero and degree one
components
$A_0$ and $A_1$, respectively,
yields a functor
from Gerstenhaber algebras to Lie-Rinehart algebras,
and that (ii) this functor has a left adjoint which assigns
the graded exterior $A$-algebra
$\Lambda_AL$ 
over $L$
to the Lie-Rinehart algebra $(A,L)$,
together with the bracket operation
(1.1)
on
$\Lambda_AL$.
Thus, given
any
Gerstenhaber algebra
$\Cal A$,
there is a canonical
morphism
of
Gerstenhaber algebras
from 
$\Lambda_{A_0}A_1$ 
to $\Cal A$.
In particular,
given the underlying $A$-module 
of $L$ which we write $L_0$ (for the moment), 
there is a bijective correspondence
between $(R,A)$-Lie algebra structures on $L_0$
and Gerstenhaber algebra structures
on
$\Lambda_AL_0$.
\smallskip
For
a general Gerstenhaber algebra $\Cal A$
over $R$,
with bracket operation written $[\cdot,\cdot]$,
an $R$-linear operator $D$
on $\Cal A$  
of degree $-1$
is said to {\it generate\/}
the Gerstenhaber bracket
provided, for every homogeneous $a, b \in \Cal A$,
$$
[a,b] = (-1)^{|a|}\left(D(ab) -(Da) b - (-1)^{|a|} a (Db)\right);
\tag1.2
$$
the operator $D$ is then called a {\it generator\/}.
For the special case where
$\Cal A = \Lambda_AL$,
$A$ and $L$
being the algebra of smooth functions and
Lie algebra of smooth vector fields on a smooth manifold,
respectively,
this terminology goes back at least to Koszul \cite\koszulon.
A general
Gerstenhaber algebra $\Cal A$
is said to be {\it exact\/}
if it is generated by an operator $D$ of square zero,
and the generator is then likewise said to be {\it exact\/};
it is then manifestly a differential
on $\Cal A$.
An exact generator will henceforth be written $\partial$.
For an exact structure,
up to a sign, the bracket operation in $\Cal A$,
viewed as an element of the usual Hom-complex
$\roman{Hom}(\Cal A \otimes_R \Cal A, \Cal A)$
with its  differential induced by that on
$\Cal A$ and that on the tensor square,
is the boundary of the multiplication map of $\Cal A$,
viewed as a chain in
$\roman{Hom}(\Cal A \otimes_R \Cal A, \Cal A)$.
A Gerstenhaber algebra with an exact generator
is called a {\it Batalin-Vilkovisky\/} algebra
\cite{\kosmathr,\,\liazutwo,\,\stashnin,\,\xuone};
in the literature,
the terminology
{\it exact\/} Gerstenhaber
algebra
occurs as well.
\smallskip
Let $M$ be an $A$-module.
Recall that a left $L$-module structure
$L\otimes_R M \to M$
on $M$, written
$(\alpha,x) \mapsto \alpha(x)$,
is called a {\it left\/}
$(A,L)$-module structure provided
$$
\align
\alpha (ax) &= \alpha(a) x + a \alpha(x),
\tag1.3.1
\\
(a\alpha)(x) &= a(\alpha (x)),
\tag1.3.2
\endalign
$$
where $a \in A,\ x \in M, \ \alpha \in L$.
More generally,
such an assignment
$L\otimes_R M \to M$,
not necessarily a left $L$-module structure
but still satisfying (1.2.1) and (1.2.2),
is referred to as an
$(A,L)$-{\it connection\/}, cf. \cite\poiscoho;
in this language,
a left $(A,L)$-module structure
is called a 
{\it flat\/} $(A,L)$-{\it connection\/}.
See \cite\poiscoho\ (2.16)
for historical remarks on these algebraic notions of connection etc.
Likewise, let $N$ be an $A$-module,
and let there be given an assignment
$N\otimes_R L \to N$,
written
$(x,\alpha) \mapsto x \circ \alpha$
or, somewhat simpler,
$(x,\alpha) \mapsto x \alpha$
(when there is no risk of confusion);
it is
called a {\it right\/}
$(A,L)$-module structure provided
it is  a right 
$L$-module structure
and, moreover, satifies
$$
\align
(ax)\alpha &= a(x \alpha) - (\alpha(a))x,
\tag1.4.1
\\
x(a\alpha) &= a(x \alpha) - (\alpha(a))x,
\tag1.4.2
\endalign
$$
where $a \in A,\ x \in N, \ \alpha \in L$;
we shall refer to such an assigment as a {\it right\/}
$(A,L)$-{\it connection\/} provided 
it only satisfies (1.3.1) and (1.3.2)
without
necessarily being a structure
of a right $L$-module.
Again, 
a right $(A,L)$-module structure
is also said to be a 
{\it flat right\/} $(A,L)$-{\it connection\/}.
\smallskip
We can now spell
our first result.

\proclaim{Theorem 1}
There is a bijective
correspondence between
right
$(A,L)$-connections
on $A$
and $R$-linear operators $D$ generating the Gerstenhaber 
bracket on 
$\Lambda_A L$.
Under this correspondence, flat 
right $(A,L)$-connections
on $A$,
that is,
right
$(A,L)$-module structures
on $A$,
correspond to operators
of square zero.
More precisely:
Given an $R$-linear operator $D$ generating the Gerstenhaber 
bracket on 
$\Lambda_A L$,
the formula
$$
a \circ \alpha = a(D\alpha) - \alpha(a),\quad
a \in A,\ \alpha \in L,
\tag1.5
$$
defines
a right
$(A,L)$-connection on $A$.
Conversely, given a
structure of 
a right $(A,L)$-connection on $A$,
the operator
$D$ on $\Lambda_AL$ defined by means of
$$
\aligned
D\langle \alpha_1, \dots ,\alpha_n \rangle
&=
\quad\sum_{i=1}^n (-1)^{(i-1)}(1\circ \alpha_i)
\langle\alpha_1, \dots\widehat{\alpha_i}\dots, \alpha_n \rangle
\\
&\phantom{=}+\quad
\sum_{j<k} (-1)^{(j+k)}
\langle \lbrack \alpha_j,\alpha_k \rbrack,
\alpha_1, \dots\widehat{\alpha_j}\dots
\widehat{\alpha_k}\dots \alpha_n\rangle.
\endaligned
\tag1.6
$$
yields an $R$-linear operator $D$ generating the Gerstenhaber 
bracket on 
$\Lambda_A L$.
\endproclaim

Before proving this Theorem we recall that
any $(R,A)$-Lie algebra
admits a universal algebra
$U(A,L)$ \cite\rinehone\ 
(the algebra of differential operators
when $A$ is the ring of smooth functions
and $L$ the Lie algebra of smooth vector fields
on a smooth manifold).
We note that left- and right $(A,L)$-modules
 plainly correspond to
left- and right $U(A,L)$-modules,
and vice versa; see 
e.~g. \cite\poiscoho\ for details.
More generally,
left- and right $(A,L)$-connections may be shown to correspond
bijectively to
left- and right $U(A,E)$-module structures,
for suitable $(R,A)$-Lie algebras $E$
mapping surjectively onto $L$.
The 
action of $L$ on $A$ which is part of the structure
of a Lie-Rinehart algebra
endows $A$ with the structure of a left
$(A,L)$- and hence
with that of a left
$U(A,L)$-module.

\demo{Proof of the first half of Theorem {\rm 1}}
Suppose given 
an $R$-linear operator $D$ generating the Gerstenhaber 
bracket on 
$\Lambda_A L$.
A straightforward verification shows that (1.5) then yields
the structure of a right $(A,L)$-connection on $A$.
Indeed, let
$x \in A,\  a \in A$ and $\alpha \in L$.
Then
$$
\align
(ax) \circ  \alpha &= (ax) (D \alpha) -\alpha(ax)
\\
&= a (x(D \alpha)) -\alpha(a)x- a(\alpha(x))
\\
&= a (x(D \alpha) - \alpha(x)) - (\alpha(a))x
\\
&= a (x \circ \alpha) - \alpha(a) x
\endalign
$$
whence (1.3.1) holds;
here
$\alpha(a)$ etc. refers to the result of
acting on $a \in A$ with $\alpha \in L$
via the left $L$-structure on $A$
(which is part of the Lie-Rinehart structure of $(A,L)$).
Moreover,
from (1.1) we deduce
$$
D(a\alpha) = a (D\alpha) + [a,\alpha]
= a (D\alpha) - \alpha(a).
$$
Hence
$$
\align
x \circ (a \alpha) &= x (D (a \alpha)) -(a \alpha)(x)
\\
&= x (a (D\alpha) - \alpha(a)) -a (\alpha(x))
\\
&= a (x (D\alpha) - \alpha(x)) -x (\alpha(a))
\\
&= a (x \circ \alpha) - (\alpha(a)) x
\endalign
$$
whence (1.3.2) holds as well.
Hence (1.5)
yields the structure of a right $(A,L)$-connection on $A$.
Another straightforward calculation
shows that the vanishing of $DD$ is equivalent to
this right connection  being flat,
that is,
to being a right $(A,L)$-module structure on $A$.
\enddemo
\smallskip
Rather than giving a direct proof for the converse,
we 
shall place the
argument to be offered in its proper context,
in the following way.
\smallskip
Recall briefly the 
Rinehart complex for $(A,L)$:
Consider the graded
left $U(A,L)$-module
$U(A,L) \otimes _A \Lambda_AL$
where $A$ acts on  $U(A,L)$ from
the right by means of the canonical map 
from $A$ to $U(A,L)$.
For $u \in U(A,L)$ and $\alpha_1, \dots, \alpha_n \in L$, let
$$
\aligned
d(u\otimes \langle \alpha_1, \dots, \alpha_n \rangle)
&=
\quad\sum_{i=1}^n (-1)^{(i-1)}u\alpha_i\otimes 
\langle\alpha_1, \dots\widehat{\alpha_i}\dots, \alpha_n \rangle
\\
&\phantom{=}+\quad
\sum_{j<k} (-1)^{(j+k)}u\otimes 
\langle \lbrack \alpha_j,\alpha_k \rbrack,
\alpha_1, \dots\widehat{\alpha_j}\dots
\widehat{\alpha_k}\dots \alpha_n\rangle.
\endaligned
\tag1.7
$$
Rinehart \cite\rinehone\ has proved that
this yields
an $R$-linear operator
$$
d
\colon
U(A,L) \otimes _A \Lambda_AL 
\longrightarrow
U(A,L) \otimes _A \Lambda_AL.
$$
We note that the non-trivial fact to be verified here is that,
for every $u \in U(A,L)$,\ $a \in A$, and $\alpha_1, \dots, \alpha_n \in L$, 
$$
d(ua\otimes \langle \alpha_1, \dots, \alpha_n \rangle)
=
d(u\otimes \langle a\alpha_1, \dots, \alpha_n \rangle)
=
\dots
=
d(u\otimes \langle \alpha_1, \dots, a\alpha_n \rangle).
\tag1.8
$$
The resulting graded object
$$
K(A,L) = (U(A,L) \otimes _A \Lambda_A(sL),d)
$$
is the
{\it Rinehart complex\/} for $(A,L)$.
Rinehart \cite\rinehone\ has also proved  that
$dd =0$, that is, $d$ is an ${U(A,L)}$-linear differential whence
$K(A,L)$ is indeed a 
chain complex, cf. also \cite\poiscoho.
When $L$ is projective as an $A$-module,
$K(A,L)$ is in fact a projective resolution of
$A$ in the category of left $U(A,L)$-modules.

\demo{Proof of the second half of Theorem 1}
Suppose given a right $(A,L)$-connection on $A$,
written $(a,\alpha) \mapsto a \circ \alpha$.
Formally the same argument
as that given by Rinehart in
\cite\rinehone\ 
establishing (1.8)
shows that
the operator
$D$ on $\Lambda_AL$ given by 
(1.6) is well defined,
that is to say,
for every $a \in A$ and $\alpha_1, \dots, \alpha_n \in L$, 
$$
D\langle a\alpha_1, \dots, \alpha_n \rangle
=
\dots
=
D\langle \alpha_1, \dots, a\alpha_n \rangle;
\tag1.9
$$
this relies on the defining properties
(1.4.1) and (1.4.2) of a right
$(A,L)$-connection.
Alternatively, when the right
$(A,L)$-connection on $A$
is flat, that is,
a right
$(A,L)$-module and hence
$U(A,L)$-module structure on $A$,
comparing
(1.6) with (1.7), 
we see that,
ignoring the algebra structure of
$\Lambda_AL$,
the operator
(1.6)
then
coincides with that of the chain complex
$(A\otimes_{U(A,L)}K(A,L),d)$
arising from taking the tensor product
of $A$ with $K(A,L)$ over $U(A,L)$,
with reference to the right
$(A,L)$-module structure on $A$.
For a general right
$(A,L)$-connection on $A$,
this kind of argument can still be used,
with $L$ being replaced by a
suitable $(R,A)$-Lie algebra $E$
mapping surjectively onto $L$,
so that the right
$(A,L)$-connection on $A$
comes from an honest right
$(A,E)$-module structure on $A$.
\smallskip
To complete the proof, we 
consider the exterior algebra $\Lambda_RL$; 
the formula (1.6) yields 
an operator 
on $\Lambda_RL$, too, which
we denote by $D$
as well, with an abuse of notation.
We then write this operator
as a sum
$$
D = D^{\wedge} + D^{[\cdot,\cdot]}
\colon
\Lambda_RL
@>>>
\Lambda_RL
$$
where,
given $\alpha_1, \dots, \alpha_n \in L$,
$$
\aligned
D^{\wedge}\langle \alpha_1, \dots, \alpha_n \rangle
&=
\quad\sum_{i=1}^n (-1)^{(i-1)}(1\circ \alpha_i)
\langle\alpha_1, \dots\widehat{\alpha_i}\dots, \alpha_n \rangle
\\
D^{[\cdot,\cdot]}
&=\sum_{j<k} (-1)^{(j+k)}
\langle \lbrack \alpha_j,\alpha_k \rbrack,
\alpha_1, \dots\widehat{\alpha_j}\dots
\widehat{\alpha_k}\dots \alpha_n\rangle.
\endaligned
$$
We note that neither
$D^{\wedge}$ nor  $D^{[\cdot,\cdot]}$
descends 
from
$\Lambda_RL$
to an operator
on
$\Lambda_AL$
but their sum does.
Likewise we write
$$
[\cdot,\cdot]
\colon
\Lambda_RL
\otimes_R
\Lambda_RL
@>>>
\Lambda_RL
$$
for the corresponding bracket.
A straightforward calculation shows that,
given $\alpha_1, \dots, \alpha_n \in L$
and
letting
$u = \langle\alpha_1, \dots, \alpha_\ell\rangle \in \Lambda^{\ell}_RL$ and
$v = \langle \alpha_{\ell+ 1}, \dots, \alpha_n\rangle \in \Lambda^{n-\ell}_RL$,
we get
$$
\aligned
D^{\wedge} (uv)
&=
(D^{\wedge}u) v
+ (-1)^{|u|} u(D^{\wedge}v) \in \Lambda_RL,
\\
D^{[\cdot,\cdot]} (uv)
&=
(D^{[\cdot,\cdot]}u) v
+ (-1)^{|u|} u(D^{[\cdot,\cdot]}v) 
+ (-1)^{|u|}[u,v] \in \Lambda_RL.
\endaligned
$$
Since the sum
$D^{\wedge} +D^{[\cdot,\cdot]}$
descends to an operator
on
$\Lambda_AL$
we conclude that,
on $\Lambda_AL$,
$$
D(uv) = (Du) v + (-1)^{|u|}u (Dv) + (-1)^{|u|} [u,v]
$$
or, equivalently,
$$
[u,v] = (-1)^{|u|} \left(D(uv) - (Du) v - (-1)^{|u|}u (Dv)\right) 
$$
as asserted. \qed
\enddemo

\smallskip
\noindent
{\smc Remark 1.}
Given a structure of a right $(A,L)$-connection
$(a,\alpha) \mapsto a \circ \alpha$
on $A$,
the square
$DD$ 
of the operator $D$
given by (1.6)
will in general be non-zero---in fact
the failure from being zero is precisely measured by 
the appropriate notion
of curvature
for this 
right $(A,L)$-connection.
We refrain from spelling out details here.

\smallskip
\noindent
{\smc Remark 2.}
For an ordinary Lie algebra
$\fra g$ over the ground ring $R$,
viewed as an $(R,R)$-Lie algebra
with $A=R$, the statement of Theorem 1 comes down 
to the  observation, cf. \cite \kosmathr,
that the complex $(\Lambda_R \fra g,\partial)$
which computes the  homology
of $\fra g$ with trivial coefficients
may be viewed
as that underlying the Batalin-Vilkovisky algebra
$(\Lambda_R \fra g,[\cdot,\cdot],\partial)$
with bracket
$[\cdot,\cdot]$
defined by (1.1).
\smallskip
Given a
generating operator $\partial$ of square zero for the Gerstenhaber 
bracket,
we shall accordingly write
$A_\partial$ for $A$ together with the 
right
$(A,L)$-module structure
given by (1.5).
An exact Gerstenhaber algebra or, equivalently,
Batalin-Vilkovisky algebra structure
on $\Lambda_A L$
has as underlying structure that of a chain complex,
and we now spell out
explicitly an observation which was already implicit
in the proof of Theorem 1.

\proclaim {Theorem 2}
Given an exact generator 
$\partial$
for the Gerstenhaber
algebra $\Lambda_AL$, 
the Batalin-Vilkovisky algebra
$(\Lambda_AL,\partial)$
coincides as a chain complex with 
\linebreak
$(A_{\partial}\otimes_{U(A,L)}K(A,L),d)$.
In particular,
when $L$ is projective as an $A$-module,
the Batalin-Vilkovisky algebra
$(\Lambda_AL,\partial)$
computes
$$
\roman H_*(L,A_\partial) \left(= \roman{Tor}_*^{U(A,L)}(A_\partial,A)\right),
$$
the homology of $L$ with coefficients in $A_\partial$. \qed
\endproclaim

In general, that is, when $L$ is not necessarily projective,
$(\Lambda_AL,\partial)$
computes a certain relative Tor-functor; see
\cite\hochsone\
for details on relative homological algebra.

\beginsection 2. Duality

We now suppose that,
as an $A$-module,
$L$ is  projective  
of finite rank (say) $n$,
so that
$\Lambda_A^nL$ is the highest non-zero exterior power of $L$
in the category of $A$-modules.
In an earlier paper \cite\duality,
we introduced
the concept of 
{\it dualizing module\/} $C_L$
for such a Lie-Rinehart algebra.
Recall that, by definition,
$C_L=\roman H^n(L, U(A,L))$, with its induced structure of a
right $(A,L)$-module (from the obvious right $U(A,L)$-module
structure on $U(A,L)$).
Moreover,
the negative of the Lie-derivative
endows
$\roman{Hom}(\Lambda_A^n L,A)$
with the structure of
a right $(A,L)$-module and,
by Theorem 2.8 of \cite\duality,
$C_L$ and  $\roman{Hom}(\Lambda_A^n L,A)$
are isomorphic as right $(A,L)$-modules.

\proclaim{Theorem 3}
There is a bijective
correspondence between
$(A,L)$-connections on $\Lambda_A^nL$
and
right
$(A,L)$-connections
on $A$.
Under this correspondence,
left $(A,L)$-module structures on $\Lambda_A^nL$
(i.~e. flat connections)
correspond to
right
$(A,L)$-module structures
on $A$.
More precisely: Given an
$(A,L)$-connection $\nabla$ on $\Lambda_A^nL$,
the negative of the (generalized) Lie-derivative
on
$A \cong \roman {Hom}_A(\Lambda_A^n L, M)
$
with reference to the
connection $\nabla$ on $M=\Lambda_A^n L$, that is, the formula
$$
(\phi \alpha) x = 
\phi(\alpha x) -\nabla_{\alpha} (\phi(x)),
\quad x \in \Lambda_A^n L,\ \alpha \in L,
\ \phi \in \roman {Hom}_A(\Lambda_A^n L, \Lambda_A^n L) \cong A,
\tag2.1
$$
yields 
the structure of 
a right $(A,L)$-connection on $A$.
Conversely,
given a
structure of 
a right $(A,L)$-connection on $A$
(written
$(a,\alpha) \mapsto a \alpha$),
on $\Lambda_A^nL \cong \roman {Hom}_A(C_L, A)$,
the assignment
$$
(\nabla_\alpha \psi) x
= \psi (x\alpha)
- (\psi x) \alpha,
\quad
x \in C_L,\ \alpha \in L,
\ \psi \in \roman {Hom}_A(C_L, A),
\tag2.2
$$
yields the structure
of an $(A,L)$-connection $\nabla$
(written
$L \otimes _R \Lambda_A^nL \to \Lambda_A^nL$,
\linebreak $(\alpha, \psi) \mapsto \nabla_\alpha \psi$).
\endproclaim

\demo{Proof}
This is straightforward and left to the reader. \qed
\enddemo

Combining Theorem 1 with Theorem 3, we obtain.

\proclaim{Corollary}
There is a bijective correspondence between
$(A,L)$-connections
on
$\Lambda_A^n L$
and linear operators $D$ generating the Gerstenhaber 
bracket on 
$\Lambda_A L$.
Under this correspondence, flat connections correspond to operators
of square zero, that is, to differentials.
The relationship is made explicit
by means of
{\rm (1.5), (1.6), (2.1)} and
{\rm (2.2)}. \qed
\endproclaim

\smallskip
\noindent
{\smc Remark.}
The statement of the corollary was proved by
Koszul \cite\koszulon\ 
(Section 2)
for the special case where
$A$ is the ring of smooth functions 
and $L$ the $(\Bobb R,A)$-Lie 
algebra of smooth vector fields on a smooth manifold,
and by Xu \cite\xuone\ (3.8)
for the special case where
$A$ is the ring of smooth functions
and $L$ the $(\Bobb R,A)$-Lie algebra of sections
of a general Lie algebroid on a smooth manifold.
However, our approach in terms of
right $(A,L)$-module structures on $A$
(coming  into play in
Theorem 1 and Theorem 3)
seems to be new
and is far more general.
\smallskip
Next we recall from
Theorem 2.8 of \cite\duality\ that
$(A,L)$ satisfies duality and inverse duality in dimension $n$;
this means that
there are  natural isomorphisms
$$
\phi \colon
\roman H^k(L,M)
@>>>
\roman H_{n-k}(L,C_L \otimes_A M)
$$
for all non-negative integers $k$ and all
left $(A,L)$-modules $M$ and, furthermore,
natural isomorphisms
$$
\psi
\colon
\roman H_{k}(L,N)
@>>>
\roman H^{n-k}(L,\roman{Hom}_A(C_L,N))
$$
for all non-negative integers $k$ and all
right $(A,L)$-modules $N$.
We now
combine this
with the results of the present paper,
in the following way:
Let
$\nabla$ be
a flat 
$(A,L)$-connection
on
$\Lambda_A^n L$; from Theorem 3 we know that
it determines
a unique
right
$(A,L)$-module 
$A_\nabla$
having $A$ as underlying $A$-module.
Further, by Theorem 1,
this right
$(A,L)$-module structure
determines a corresponding exact generator
$\partial$
for the Gerstenhaber algebra $\Lambda_A L$ and,
by Theorem 2,
as a chain complex,
the resulting
Batalin-Vilkovisky algebra
$(\Lambda_A L,\partial)$
coincides
with
the chain complex $(A_\partial \otimes_{U(A,L)}K(A,L),d)$
arising from the Rinehart resolution 
$K(A,L)$ of $A$ in the category of left $(A,L)$-modules and
computing
$$
\roman H_*(L,A_\nabla)
\quad (=\roman{Tor}_*^{U(A,L)}(A_\nabla,A)),
$$
the homology of $L$
with coefficients in the right
$(A,L)$-module 
$A_\nabla$.

\proclaim {Theorem 4}
The homology of $L$
with coefficients in the right
$(A,L)$-module 
$A_\nabla$,
that is,
that of the Batalin-Vilkovisky algebra
$(\Lambda_A L,\partial)$,
is naturally isomorphic to
the cohomology
$$
\roman H^*(L,\Lambda_A^n L_\nabla)
\quad \left(=\roman{Ext}^*_{U(A,L)}(A,\Lambda_A^n L_\nabla)\right)
$$
of 
$L$ with coefficients in $\Lambda_A^n L_\nabla$,
the $A$-module
$\Lambda_A^n L$,
with left $(A,L)$-module structure
given by the connection $\nabla$. \qed
\endproclaim

\beginsection 3. Some applications and concluding remarks

(3.1) {\smc Lie algebroids.} 
Let $(A,L)$ be the Lie-Rinehart algebra
arising from a Lie algebroid over a smooth manifold
of dimension $n$,
so that
$A$ is the algebra of smooth functions
and the $L$ underlying $A$-module the space of sections of the
Lie algebroid.
Then 
$n$ is the rank of $L$, and
$\Lambda_A^n L$,
being the space of sections of a real line bundle,
manifestly has a flat $(A,L)$-connection $\nabla$.
(Over the complex numbers, in particular in the holomorphic context,
there will in general be obstructions to the existence
of such a flat connection.)
Hence $A$ has the structure of a 
right
$(A,L)$-module (Theorem 3),
to be referred to as $A_\nabla$
which, by Theorem 1, 
determines a generator $\partial$
for the Gerstenhaber algebra
$\Lambda_AL$; 
by Theorem 2, the homology
of the resulting Batalin-Vilkovisky algebra
$(\Lambda_AL,\partial)$ 
is naturally isomorphic
to the 
homology of $L$ with coefficients
in $A_\nabla$ and,
by Theorem 4, to the
cohomology
of $L$ with values in
$\Lambda_A^n L_\nabla$,
the $A$-module
$\Lambda_A^n L$,
with left $(A,L)$-module structure
given by the connection $\nabla$.
The isomorphism between
the homology  of
$(\Lambda_AL,\partial)$ 
and the cohomology
of $L$ with values in
$\Lambda_A^n L_\nabla$
has been observed in
\cite\xuone\ (4.6)
for the special case
where
the line bundle
underlying 
$\Lambda_A^n L$
is trivial 
in such a way that,
as left $(A,L)$-modules,
$\Lambda_A^n L$
and $A$ are isomorphic or,
more formally, 
when the modular class of
the Lie-Rinehart algebra $(A,L)$ is trivial,
cf.
\cite \duality\
and also
\cite{\evluwein, \weinsfte}.
Xu's isomorphism (4.6)
is in fact a special case
of the duality isomorphism
$$
\roman H_*(L,A_\nabla) \cong \roman H^{n-*}(L, \Lambda^n_A L_\nabla)
\tag3.1.1
$$
established in our paper \cite\duality.
Notice that when $L$ is the $(\Bobb R,A)$-Lie algebra
of smooth vector fields,
we can view the elements of
$\Lambda_AL$ as de Rham currents
and hence
the Batalin-Vilkovisky algebra
$(\Lambda_AL,\partial)$
as 
a chain complex 
defining 
{\it de Rham homology\/}.

\smallskip\noindent
(3.2) {\smc Poisson structures.} 
Let $A$ be a Poisson algebra, with Poisson structure
$\{\cdot,\cdot\}$,
and let $D_{\{\cdot,\cdot\}}$
be its 
$(R,A)$-Lie algebra,
coming from the Poisson structure,
see \cite\poiscoho\ (3.8) for details and also
\cite\duality.
Besides its structure
of a left
$(A,D_{\{\cdot,\cdot\}})$-module,
from the Poisson structure,
the algebra $A$ inherits also structure of a right
$(A,D_{\{\cdot,\cdot\}})$-module
by means of the formula
$$
a(b(du)) = \{ab,u\}
$$
which 
we refer to by the notation
$A_{\{\cdot,\cdot\}}$;
here $a, b, u \in A$ and $du$ denotes the differential of $u$.
By Theorem 1, this structure induces an exact generator
 for the 
Gerstenhaber algebra
$\Lambda_A D_{\{\cdot,\cdot\}}$;
we denote
this generator by
$\partial_{\{\cdot,\cdot\}}$.
The resulting Batalin-Vilkovisky algebra
$(\Lambda_A D_{\{\cdot,\cdot\}},\partial_{\{\cdot,\cdot\}})$
computes the Poisson
homology of $A$,
cf. \cite\poiscoho.
This follows at once from Theorem 2.
In fact
the chain complex
which underlies 
$(\Lambda_A D_{\{\cdot,\cdot\}},\partial_{\{\cdot,\cdot\}})$
coincides with
the chain complex
$(A_{\{\cdot,\cdot\}}
\otimes_{U(A,D_{\{\cdot,\cdot\}})}K(A,D_{\{\cdot,\cdot\}}),d)$
and, by definition,
the Poisson homology of $A$ is the homology of this complex
\cite\poiscoho.
For the special case
where, 
as an $A$-module,
$D_{\{\cdot,\cdot\}}$
is finitely generated projective of rank $n$,
the corresponding left
$(A,D_{\{\cdot,\cdot\}})$-module structure
on the top exterior power
$\Lambda_A^nD_{\{\cdot,\cdot\}}$
has been introduced
in (5.8) of our paper \cite\duality.
When $A$ is the Poisson algebra of smooth functions on a 
smooth Poisson manifold
and 
$D_{\{\cdot,\cdot\}}$
the corresponding $(\Bobb R,A)$-Lie algebra
arising from the cotangent Lie algebroid,
the chain complex
$(A_{\{\cdot,\cdot\}}
\otimes_{U(A,D_{\{\cdot,\cdot\}})}K(A,D_{\{\cdot,\cdot\}}),d)$
coincides with Koszul's
complex defining Poisson homology \cite\koszulon.
In fact, Koszul 
has indeed shown that then the operator
$\partial_{\{\cdot,\cdot\}}$
is a generator for the corresponding
Gerstenhaber algebra.
For this special case, a geometric description
of the corresponding left $(A,D_{\{\cdot,\cdot\}})$-module
structure 
on 
$\Lambda^n_AD_{\{\cdot,\cdot\}}$
or, rather, that of the corresponding flat connection
on the underlying line bundle
was given in
\cite\evluwein;
an expression relating
the latter
with the generator of
the corresponding
Batalin-Vilkovisky algebra
has been given in (18) of \cite\xuone. 
In (4.8) of that paper,
a certain isomorphism between Poisson
homology and cohomology is given which
is again a special case of (3.1.1)
above;
see our paper
\cite\duality\ (5.4) for details.
\smallskip
Let $L = \roman{Der}(A)$, with its obvious
structure of an $(R,A)$-Lie algebra.
The pair $(L,D_{\{\cdot,\cdot\}})$
constitutes a structure of what we shall call
elsewhere an
$(R,A)$-{\it Lie bialgebra\/},
in fact,
a {\it triangular\/}
$(R,A)$-{\it  Lie bialgebra\/}.
These are abstractions from the notion of a
Lie bialgebroid
and from that of a triangular
one,  respectively,
cf. \cite\kosmathr.
We only note here
that the notion of an 
$(R,A)$-Lie bialgebra
is equivalent to that of a strong differential
Gerstenhaber algebra,
that is, to a Gerstenhaber algebra
endowed with a derivation of degree 1 and square
zero which behaves as a derivation
for the algebra and bracket structures;
cf. \cite\kosmathr\  and \cite\xuone.
\smallskip\noindent
(3.3) Finally we comment on question 2 
and question 3 in Section 5 of \cite\xuone:
\newline\noindent
(3.3.1) Given a general Lie algebroid,
for the corresponding Gerstenhaber algebra,
there does not seem to exist
a canonical generating operator 
corresponding to its modular class:
We know from Theorem 1 above that,
writing $(A,L)$ for the corresponding Lie-Rinehart algebra,
these generating operators correspond
bijectively to right $(A,L)$-module structures
(or connections)
on $A$
and, for a Poisson algebra $A$, the Poisson
structure
(somewhat miracously)
induces a right module structure
for the corresponding
Lie-Rinehart algebra as well
and hence a canonical exact generator.
Now, on the one hand, by Theorem  3,
right $(A,L)$-module structures
on $A$
correspond bijectively
to left
$(A,L)$-module structures
on the top exterior power $\Lambda_A^nL$ of $L$.
On the other hand, the modular class
is determined by
the 
left
$(A,L)$-module
$Q_L = \roman {Hom}_A(C_L, \omega_A)$
(see \cite\duality\ for the 
right $(A,L)$-module  $\omega_A$),
but there is no obvious way
to relate
these structures with requisite
right
$(A,L)$-module
structures 
(or right connections)
on $A$.
The fact that there is a canonical such structure
in the Poisson case seems to reflect certain compatibility
properties between left and right
structures in this case.
Indeed, the Koszul operator 
$\partial_{\{\cdot,\cdot\}}$
on $\Lambda_A D_{\{\cdot,\cdot\}}$
looks canonical but depends
tacitly also on the right module structure
on $A$ which, in turn, depends on the Poisson structure, too.
Hence the answer to Xu's question 2
will presumably be No.
\newline\noindent
(3.3.2) The family of Poisson homologies
parametrized by the first Poisson cohomology
are precisely the homologies
$\roman H_*^{\roman{Poisson}}(A,A_\nabla)$
where $\nabla$ runs
through
connections
on $\Lambda_A^n L$;
these homologies
are isomorphic to
the cohomologies
$\roman H^*_{\roman{Poisson}}(A,\Lambda_A^nL_\nabla)$,
where
$L= D_{\{\cdot,\cdot\}}$,
the corresponding 
$(\Bobb R,A)$-Lie algebra
having rank $n$ as a projective $A$-module.
This answers  Xu's question 3.
\bigskip
\widestnumber\key{999}
\centerline{References}
\smallskip\noindent

\ref \no \batviltw
\by I. A. Batalin and G. S. Vilkovisky
\paper Quantization of gauge theories
with linearly dependent generators
\jour  Phys. Rev. 
\vol D 28
\yr 1983
\pages  2567--2582
\endref

\ref \no \batvilfo
\by I. A. Batalin and G. S. Vilkovisky
\paper Closure of the gauge algebra, generalized Lie equations
and Feynman rules
\jour  Nucl. Phys. B
\vol 234
\yr 1984
\pages  106-124
\endref

\ref \no \batavilk
\by I. A. Batalin and G. S. Vilkovisky
\paper Existence theorem for gauge algebra
\jour Jour. Math. Phys.
\vol 26
\yr 1985
\pages  172--184
\endref

\ref \no \evluwein
\by S. Evens, J.-H. Lu, and A. Weinstein
\paper Transverse measures, the modular class, and a cohomology pairing
for Lie algebroids
\paperinfo preprint
\endref

\ref \no \gersthtw
\by M. Gerstenhaber
\paper The cohomology structure of an associative ring
\jour Ann. of Math.
\vol 78
\yr 1963
\pages  267-288
\endref

\ref \no \geschthr
\by M. Gerstenhaber and Samuel D. Schack
\paper Algebras, bialgebras, quantum groups and algebraic
deformations
\paperinfo In: Deformation theory and quantum groups with
applications to mathematical physics, M. Gerstenhaber and J. Stasheff, eds.
\jour Cont. Math.
\vol 134
\pages 51--92
\publ AMS
\yr 1992
\publaddr Providence 
\endref

\ref \no \getzltwo
\by E. Getzler
\paper Batalin-Vilkovisky algebras and two-dimensional topological field
theories
\jour Comm. in Math. Phys.
\vol 195
\yr 1994
\pages 265--285
\endref

\ref \no \hochsone
\by G. Hochschild
\paper Relative homological algebra
\jour  Trans. Amer. Math. Soc.
\vol 82
\yr 1956
\pages 246--269
\endref

\ref \no \poiscoho
\by J. Huebschmann
\paper Poisson cohomology and quantization
\jour J. f\"ur die Reine und Angew. Math.
\vol 408
\yr 1990
\pages 57--113
\endref
\ref \no \duality
\by J. Huebschmann
\paper 
Duality for Lie-Rinehart algebras and the modular class
\paperinfo preprint, 1997,
dg-ga/9702008
\endref

\ref \no \husmosta
\by D. Husemoller, J. C. Moore, and J. D. Stasheff
\paper Differential homological algebra and homogeneous spaces
\jour J. of Pure and Applied Algebra
\vol 5
\yr 1974
\pages  113--185
\endref

\ref \no \kosmathr
\by Y. Kosmann-Schwarzbach 
\paper Exact Gerstenhaber algebras and Lie bialgebroids
\jour  Acta Applicandae Mathematicae
\vol 41
\yr 1995
\pages 153--165
\endref

\ref \no \koszulon
\by J. L. Koszul
\paper Crochet de Schouten-Nijenhuis et cohomologie
\jour Ast\'erisque,
\vol hors-serie,
\yr 1985
\pages 251--271
\paperinfo in E. Cartan et les Math\'ematiciens d'aujourd'hui, 
Lyon, 25--29 Juin, 1984
\endref

\ref \no \liazutwo
\by B. H. Lian and G. J. Zuckerman
\paper New perspectives on the BRST-algebraic structure
of string theory
\jour Comm. in Math. Phys.
\vol 154
\yr 1993
\pages  613--646
\endref

\ref \no \rinehone
\by G. Rinehart
\paper Differential forms for general commutative algebras
\jour  Trans. Amer. Math. Soc.
\vol 108
\yr 1963
\pages 195--222
\endref

\ref \no \stashnin
\by J. D. Stasheff
\paper Deformation theory and the Batalin-Vilkovisky 
master equation
\paperinfo in: Ascona, 1996, to appear
\endref

\ref \no  \weinsfte
\by A. Weinstein
\paper The modular automorphism group of a Poisson manifold
\paperinfo to appear in: special volume in honor of A. Lichnerowicz
\jour J. of Geometry and Physics
\endref

\ref \no \xuone
\by Ping Xu
\paper 
Gerstenhaber algebras and BV-algebras
in Poisson geometry
\paperinfo preprint, 1997
\endref

\enddocument